\documentclass[11pt]{amsart}

\usepackage[dvipsnames]{xcolor}

\usepackage[
		backend=biber,
		style=alphabetic,
		doi=true,
		isbn=false,
		url=true,
		maxbibnames=99,
		maxalphanames=6,
		minalphanames=3
	]{biblatex}
\DeclareFieldFormat[
		article,
		book,
		inbook,
		incollection,
		inproceedings,
		patent,
		thesis,
		unpublished
	]{title}{\emph{#1\isdot}}
\DeclareFieldFormat{journaltitle}{#1\isdot}
\AtEveryBibitem{
	\clearfield{note}
	\iffieldundef{doi}{}{\clearfield{url}}
	\iffieldequalstr{eprinttype}{arxiv}{\clearfield{url}}{}
	\iffieldequalstr{eprinttype}{arXiv}{\clearfield{url}}{}
}

\usepackage{amssymb, mathrsfs, mathtools, braket, mleftright}

\usepackage{tikz}
\usetikzlibrary{cd}

\usepackage[
		colorlinks,
		linkcolor=Maroon,
		citecolor=ForestGreen,
		urlcolor=MidnightBlue
	]{hyperref}

\usepackage[nameinlink]{cleveref}

\crefname{section}{Section}{Sections}
\crefname{subsection}{Section}{Sections}
\crefname{subsubsection}{Section}{Sections}

\makeatletter
\def\l@section{\@tocline{1}{0pt}{1pc}{}{\bfseries}}
\def\l@subsection{\@tocline{2}{0pt}{2.5pc}{5pc}{}}
\renewcommand{\subsubsection}{%
	\@startsection{subsubsection}{3}{\z@}%
		{.5\linespacing\@plus.7\linespacing}{-.5em}%
		{\normalfont\bfseries}}
\makeatother

\newtheorem{thm}{Theorem}[section]
\crefname{thm}{Theorem}{Theorems}
\newtheorem{prop}[thm]{Proposition}
\crefname{prop}{Proposition}{Propositions}
\newtheorem{lem}[thm]{Lemma}
\crefname{lem}{Lemma}{Lemmas}
\newtheorem{cor}[thm]{Corollary}
\crefname{cor}{Corollary}{Corollaries}

\crefname{conj}{Conjecture}{Conjectures}

\crefname{ques}{Question}{Questions}

\theoremstyle{definition}
\newtheorem{defi}[thm]{Definition}
\crefname{defi}{Definition}{Definitions}
\newtheorem{setting}[thm]{Setting}
\crefname{setting}{Setting}{Settings}
\newtheorem{exa}[thm]{Example}
\crefname{exa}{Example}{Examples}
\newtheorem{claim}{Claim}
\crefname{claim}{Claim}{Claims}

\theoremstyle{remark}
\newtheorem{rem}[thm]{Remark}
\crefname{rem}{Remark}{Remarks}

\crefname{enumi}{}{}
\creflabelformat{enumi}{#2(#1)#3}

\crefname{enumii}{}{}
\creflabelformat{enumii}{#2(#1)#3}

\numberwithin{equation}{section}

\newcommand{\abs}[1]{\lvert #1\rvert}
\newcommand{\card}{\abs}
\newcommand{\isom}{\cong}
\newcommand{\restr}[2]{#1\rvert_{#2}}
\newcommand{\bigrestr}[2]{#1\bigr\rvert_{#2}}
\newcommand{\isomarrow}{\xrightarrow{\;\sim\;}}

\newcommand{\PP}{\mathbb{P}}
\renewcommand{\AA}{\mathbb{A}}
\newcommand{\ZZ}{\mathbb{Z}}
\newcommand{\QQ}{\mathbb{Q}}
\newcommand{\RR}{\mathbb{R}}
\newcommand{\CC}{\mathbb{C}}
\newcommand{\Oc}{\mathcal{O}}
\newcommand{\Ic}{\mathcal{I}}

\newcommand{\Ec}{\mathcal{E}}
\newcommand{\Fc}{\mathcal{F}}
\newcommand{\Lc}{\mathcal{L}}

\newcommand{\Sc}{\mathcal{S}}

\DeclareMathOperator{\Bl}{Bl}

\DeclareMathOperator{\GL}{GL}
\DeclareMathOperator{\Gr}{Gr}
\DeclareMathOperator{\Hom}{Hom}
\renewcommand{\ker}{\operatorname{Ker}}

\DeclareMathOperator{\Proj}{Proj}
\DeclareMathOperator{\rank}{rank}

\DeclareMathOperator{\Sym}{Sym}

\newcommand{\sset}[2]{\set{#1|#2}}

\title[Fano blowups of Grassmannians]{Fano blowups of Grassmannians along sub-Grassmannians}
\author{Yuto Masamura}
\address{Graduate School of Mathematical Sciences, The University of Tokyo, 3-8-1 Komaba, Meguro-ku, Tokyo 153-8914, Japan}
\email{masamura@ms.u-tokyo.ac.jp}

\begin{document}

\begin{abstract}
	Let $\tilde X$ be the successive blowup of the Grassmannian $\Gr(d,n+1)$ with $n\ge3$, in arbitrary characteristic, along sub-Grassmannians of codimension $r_1,\dotsc,r_s\ge2$ that meet transversely.
	We prove that, if $d\ne2$, then $\tilde X$ is Fano if and only if $\sum_i(r_i-1)\le n$, and we give a similar criterion for $d=2$.
	As applications, we prove that the successive blowup of $\PP^4$ along $s$ planes that meet transversely is Fano if and only if $s\le4$.
	We also determine when the successive blowup of a smooth quadric fourfold along planes of one family that meet transversely is Fano.
\end{abstract}

\maketitle

\tableofcontents

\section{Introduction}\label{sec:intro}

Fano varieties are basic objects in birational geometry, and their classification in each dimension is a central problem.
In dimension four, the classification of smooth Fano varieties is open.
Casagrande studies smooth Fano fourfolds $X$ over $\CC$ with large Picard number $\rho(X)$ \cite{Casagrande}.
She proves that if $\rho(X)>9$, then $X$ is a product of del Pezzo surfaces \cite[Theorem~1.1]{Casagrande}, and she studies in detail the fourfolds with $\rho(X)\in\{7,8,9\}$, for which only a few families are known.
In particular, she proves that if $\rho(X)\ge7$ and $X$ is not a product of surfaces and has no small elementary contraction, then $X$ is obtained from a cubic fourfold by successively blowing up $\rho(X)-1$ planes that meet pairwise in a point \cite[Theorem~1.2]{Casagrande}.
Thus successive blowups of fourfolds along planes that meet pairwise in a point appear naturally in the classification.
In the last section of \cite{Casagrande}, she considers such blowups of a smooth cubic fourfold, together with the analogous blowups of $\PP^4$ and of a smooth quadric fourfold $Q$, and asks when these fourfolds are Fano.

The projective space $\PP^4$ is the Grassmannian $\Gr(1,5)$ of one-dimensional subspaces of $\CC^5$, and $Q$ is the Grassmannian $\Gr(2,4)$.
In both cases, the planes in Casagrande's questions are sub-Grassmannians, up to a duality for $Q$.
It is therefore natural to study blowups of Grassmannians along sub-Grassmannians.
The case of one center is known over $\CC$.
Hu, Ke, Li, and Song prove that the blowup of $\Gr(d,m)$ along a smooth Schubert variety, for example a sub-Grassmannian, is Fano if and only if the codimension of the center is at most $m$ \cite[Theorem~1.1]{HuKeLiSong}.
Biswas and Saha describe more generally when the blowup of a generalized flag variety along a smooth Schubert variety is Fano \cite[Corollary~3.11]{BiswasSaha}.

In this paper, we generalize the blowup of a Grassmannian along a sub-Grassmannian to an arbitrary number of centers.
We determine when the successive blowup of a Grassmannian along sub-Grassmannians that meet transversely is Fano.
As applications, we answer Casagrande's questions for $\PP^4$ and for $Q$.
The remaining question, for cubic fourfolds \cite[Question~11.6]{Casagrande}, has recently been answered by Secci, who constructs Fano blowups along up to $8$ planes \cite[Theorem~1.1]{Secci}.

Throughout this paper, we work over an algebraically closed field $k$ of arbitrary characteristic.
Let $V$ be a vector space of dimension $n+1$, let $1\le d\le n$, and let $X\coloneq\Gr(d,V)$ be the Grassmannian of $d$-dimensional subspaces of $V$.
Let $W_1,\dotsc,W_s\subseteq V$ be subspaces of dimension at least $d$, and let $c_i$ be the codimension of $W_i$ in $V$.
The sub-Grassmannians $Z_i\coloneq\Gr(d,W_i)\subseteq X$ have codimension $r_i\coloneq dc_i$, and we assume that $r_i\ge2$.
Assume that $Z_1,\dotsc,Z_s$ meet transversely, that is, for every subset $I\subseteq\{1,\dotsc,s\}$ such that $\bigcap_{i\in I}W_i$ has dimension at least $d$, the codimension of this intersection in $V$ is $\sum_{i\in I}c_i$ (\cref{lem:subgrass}).
Let $\pi\colon\tilde X\to X$ be the blowup along $Z_1$, followed by the blowups along the strict transforms of $Z_2,\dotsc,Z_s$.
The variety $\tilde X$ is smooth, has Picard number $s+1$, and does not depend on the order of the blowups (\cref{lem:basic,lem:local}).
Our main result is the following.

\begin{thm}[{\cref{cor:grass}}]\label{introthm:grass}
	Let $V$ be a vector space of dimension $n+1\ge4$, and let $1\le d\le n$.
	Let $W_1,\dotsc,W_s\subseteq V$ be subspaces of dimension at least $d$ such that the sub-Grassmannians $\Gr(d,W_i)\subseteq\Gr(d,V)$ have codimension $r_1,\dotsc,r_s\ge2$ and meet transversely.
	Let $\tilde X$ be the successive blowup of $\Gr(d,V)$ along them.
	\begin{enumerate}
		\item If $d\ne2$, then $\tilde X$ is Fano if and only if $\sum_i(r_i-1)\le n$.
		\item If $d=2$, then $\tilde X$ is Fano if and only if $\sum_i(r_i-1)\le2n+1$ and $\sum_{i\in I}(r_i-1)\le n$ for every subset $I\subseteq\{1,\dotsc,s\}$ with $\bigcap_{i\in I}W_i\ne0$.
	\end{enumerate}
\end{thm}

Both cases come from a single criterion, which says that $\tilde X$ is Fano if and only if $\theta<n+1$, where $\theta$ is a rational number determined by the subspaces $W_i$ (\cref{thm:grass}).
For a single center ($s=1$), \cref{introthm:grass} says that $\tilde X$ is Fano if and only if $r_1\le\dim V$, which agrees with \cite[Theorem~1.1]{HuKeLiSong} over $\CC$.
The assumption $n\ge3$ is needed (\cref{rem:grass-n2}).

For $d=1$, the variety $X$ is the projective space $\PP^n$, and the centers are linear subspaces of codimension $r_i$.
In this case, if $n\ge3$, then $\tilde X$ is Fano if and only if $\sum_i(r_i-1)\le n$ by \cref{introthm:grass}.
For $n=3$, the centers are points and lines, and the Fano threefolds that we obtain appear in the classification of Mori and Mukai \cite{MoriMukai}.

We now return to Casagrande's questions.
For planes in $\PP^4$, meeting transversely means that any two of them meet in a point and no three of them have a common point (\cref{ss:grass-fourfolds}).
In particular, general planes in $\PP^4$ meet transversely.
For planes $Z_1,\dotsc,Z_s\subseteq\PP^4$, we write $\tilde X_s$ for the successive blowup of $\PP^4$ along them.
For general planes in $\PP^4$, Casagrande observes that $\tilde X_1$ and $\tilde X_2$ are toric Fano fourfolds \cite[Section~11.5]{Casagrande}.
She also proves that if $\tilde X_s$ is Fano, then $s\le4$, by applying the classification of Fano threefolds to the general fiber of a morphism $\tilde X_s\to\PP^1$ \cite[Lemma~11.3]{Casagrande}.
She asks whether $\tilde X_3$ and $\tilde X_4$ are Fano \cite[Question~11.4]{Casagrande}.
We obtain the following answer.

\begin{thm}[{\cref{cor:P4}}]\label{introthm:P4}
	Let $Z_1,\dotsc,Z_s\subseteq\PP^4$ be planes that meet transversely.
	Then the successive blowup $\tilde X_s$ of $\PP^4$ along them is Fano if and only if $s\le4$.
\end{thm}

In particular, the answer to \cite[Question~11.4]{Casagrande} is positive.
Moreover, \cref{introthm:P4} replaces the generality of the planes by transversality, and holds in arbitrary characteristic.
\Cref{introthm:grass} also applies to centers of different dimensions.
For example, the successive blowup of $\PP^4$ along two planes and a line that meet transversely is Fano, with Picard number $4$ and anticanonical degree $293$ (\cref{exa:P4-line}).

The planes in the smooth quadric fourfold $Q\isom\Gr(2,4)$ form two families, parametrized by the points and by the planes of $\PP^3$.
After applying a duality of $\PP^3$, planes $Z_1,\dotsc,Z_s\subseteq Q$ of one family correspond to points $p_1,\dotsc,p_s\in\PP^3$.
The planes meet transversely if and only if the points are pairwise distinct and no three of them are collinear (\cref{ss:grass-fourfolds}).
In this case, we write $\tilde X_s$ for the successive blowup of $Q$ along them.
For general planes of one family, Manivel proves that $\tilde X_5$ is Fano \cite{Manivel}.
Casagrande asks whether $\tilde X_s$ is Fano for $s\in\{6,7,8\}$ \cite[Question~11.5]{Casagrande}.
The case $n=3$, $d=2$, and $r_i=2$ for every $i$ in \cref{introthm:grass} gives the following answer.

\begin{thm}[{\cref{cor:quadric}}]\label{introthm:quadric}
	Let $Q$ be a smooth quadric fourfold, and let $Z_1,\dotsc,Z_s\subseteq Q$ be planes of one family that meet transversely, corresponding to points $p_1,\dotsc,p_s\in\PP^3$.
	Let $\tilde X_s$ be the successive blowup of $Q$ along $Z_1,\dotsc,Z_s$.
	Then $\tilde X_s$ is Fano if and only if $s\le7$ and no four of the points $p_i$ lie on a plane.
\end{thm}

For general points, no four of them lie on a plane.
Thus the answer to \cite[Question~11.5]{Casagrande} is positive for $s=6,7$ and negative for $s=8$.
For $s=5$, \cref{introthm:quadric} agrees with the result of Manivel \cite{Manivel}.

Over $\CC$, the fourfolds $\tilde X_6$ and $\tilde X_7$ of \cref{introthm:quadric} form new families of Fano fourfolds with Picard numbers $7$ and $8$.
They have anticanonical degrees $122$ and $78$ \cite[Table~11.2]{Casagrande}, so they are not products of surfaces, since a product of del Pezzo surfaces of degrees $d_1,d_2\le9$ has anticanonical degree $6d_1d_2$.
These degrees differ from those of the other non-product families with these Picard numbers that we know of, namely the Fano models of the blowups of $\PP^4$ at $6$ and $7$ general points \cite[Example~7.2 and Table~7.1]{CasagrandeSecci}, the blowups of the Fano model at $5$ general points along the transforms of a cubic scroll and of a sextic K3 surface \cite[Proposition~1.9 and Tables~7.2 and~7.3]{CasagrandeSecci}, and the blowups of cubic fourfolds along $6$ and $7$ planes \cite[Table~7.1]{Secci}.

We now outline the proofs of these results.
Using the results of the author on successive blowups \cite{MasamuraBlowup}, we show that the successive blowup along centers that meet transversely is the fiber product over $X$ of the blowups along the single centers (\cref{lem:local}).
This allows us to lift curves on $X$ to $\tilde X$ and control their intersection numbers with the exceptional divisors (\cref{lem:lift}).
To prove the ampleness of $-K_{\tilde X}$ in \cref{introthm:grass}, we pull back the dual of the tautological subbundle on $X=\Gr(d,V)$ to (the normalization of) each curve on $\tilde X$ and compare its Harder--Narasimhan filtration with the subspaces $W_i$.
This shows that $-K_{\tilde X}-(n+1-\theta)\pi^*\Oc_X(1)$ is nef, and hence that $-K_{\tilde X}$ is ample for $\theta<n+1$.
For the converse, we lift lines meeting the centers, which we find by Schubert calculus.
For $d=2$, we also lift morphisms $\PP^1\to X$ of degree two. When $r_i=2$ for every $i$, we obtain them from quadratic forms on $V^\vee$.

\subsubsection*{Organization of the paper}

In \cref{sec:setup}, we fix notation and prove basic facts on successive blowups of a smooth projective variety along centers that meet transversely.
In \cref{sec:grass}, we apply these results to sub-Grassmannians of a Grassmannian, prove \cref{introthm:grass}, and deduce \cref{introthm:P4,introthm:quadric}.

\subsubsection*{Acknowledgements}

The author would like to thank his advisor, Keiji Oguiso, for his guidance and constant encouragement.
The author thanks Cinzia Casagrande for suggesting that he include \cref{exa:P4-line}.
The author was supported by JSPS KAKENHI Grant Number JP24KJ0856.
This research was also supported in part by access to OpenAI models provided through the ChatGPT program for academic researchers.

\subsubsection*{AI disclosure}

The author used Claude (Anthropic) and GPT (OpenAI) for help with the research and the writing of this paper.
The author has checked all the mathematics in this paper and takes full responsibility for its content.

\section{Preliminaries and successive blowups}\label{sec:setup}

In this section, we fix notation and study successive blowups of a smooth projective variety.

\subsection{Notation and conventions}\label{ss:conventions}

We work over an algebraically closed field $k$.
A \emph{variety} is an integral separated scheme of finite type over $k$.

Let $X$ be a smooth projective variety.
We say that $X$ is \emph{Fano} if $-K_X$ is ample.
We write $N^1(X)$ for the $\RR$-vector space of $\RR$-divisors on $X$ modulo numerical equivalence, and $\rho(X)\coloneq\dim_\RR N^1(X)$ for the Picard number of $X$.

For a locally free sheaf $\Ec$ of finite rank on a scheme $X$, we use Grothendieck's convention: the projective bundle
\[
	\PP_X(\Ec)\coloneq\Proj_X\Sym\Ec
\]
parametrizes invertible quotients of $\Ec$, and $\Oc_{\PP_X(\Ec)}(1)$ denotes the tautological invertible quotient.
For a vector space $V$, we write $[v]\in\PP(V^\vee)$ for the point given by a nonzero vector $v\in V$.

For a vector space $V$ and an integer $1\le d<\dim V$, we write $\Gr(d,V)$ for the Grassmannian of $d$-dimensional subspaces of $V$, and $\Sc\subseteq V\otimes\Oc_{\Gr(d,V)}$ for its tautological subbundle of rank $d$.
We write $[P]\in\Gr(d,V)$ for the point given by a $d$-dimensional subspace $P\subseteq V$.
We embed $\Gr(d,V)$ in $\PP((\bigwedge^dV)^\vee)$ by the Pl\"ucker embedding $[P]\mapsto[\bigwedge^dP]$, so that $\Oc(1)=\det\Sc^\vee$.
For a subspace $W\subseteq V$ of dimension at least $d$, we regard $\Gr(d,W)$ as a closed subvariety of $\Gr(d,V)$.

For $d=1$, we identify $\Gr(1,V)$ with $\PP(V^\vee)$ by the Pl\"ucker embedding.
Then $\Sc=\Oc(-1)$, and $\Gr(1,W)=\PP(W^\vee)$ for every nonzero subspace $W\subseteq V$.

\subsection{Successive blowups}\label{ss:blowups}

In this subsection, we collect basic facts on successive blowups along centers that meet transversely.

We first define the transversality condition used in this paper.

\begin{defi}\label{defi:transverse}
	Let $X$ be a smooth variety, and let $Z_1,\dotsc,Z_s\subsetneq X$ be smooth subvarieties.
	The subvarieties $Z_1,\dotsc,Z_s$ \emph{meet transversely} if, at every closed point $p\in\bigcup_iZ_i$, there is a regular system of parameters $x_1,\dotsc,x_n$ of $\Oc_{X,p}$ and there are pairwise disjoint subsets $J_i\subseteq\{1,\dotsc,n\}$ such that every $Z_i$ containing $p$ satisfies
	\[
		\Ic_{Z_i,p}=(x_j\mid j\in J_i).
	\]
\end{defi}

Note that if $Z_1,\dotsc,Z_s$ meet transversely, then no $Z_i$ is contained in another $Z_h$.
This allows us to blow up their strict transforms one after another, and we fix the following setting.

\begin{setting}\label{set:main}
	Let $X$ be a smooth projective variety, and let $Z_1,\dotsc,Z_s\subseteq X$ be subvarieties of codimension $r_1,\dotsc,r_s\ge2$ that meet transversely.
	Let
	\[
	\begin{tikzcd}
		\tilde X\coloneq X^{(s)}\arrow[r,"\pi_s"]
			& X^{(s-1)}\arrow[r]
			& \dotsb\arrow[r]
			& X^{(1)}\arrow[r,"\pi_1"]
			& X^{(0)}\coloneq X
	\end{tikzcd}
	\]
	be the sequence of blowups defined by
	\[
		X^{(i)}\coloneq\Bl_{Z_i^{(i-1)}}X^{(i-1)}\qquad(1\le i\le s),
	\]
	where $Z_i^{(i-1)}\subseteq X^{(i-1)}$ denotes the strict transform of $Z_i$.
	Set $\pi\coloneq\pi_1\circ\dotsb\circ\pi_s\colon\tilde X\to X$, and let $E_i$ be the pullback to $\tilde X$ of the exceptional divisor $E_i^{(i)}\subseteq X^{(i)}$ of $\pi_i$.
	For a closed point $p\in X$, set $I(p)\coloneq\sset{i}{p\in Z_i}$.
	For each $i$, let $\pi'_i\colon\Bl_{Z_i}X\to X$ be the blowup along $Z_i$, and let $\bar E_i\subseteq\Bl_{Z_i}X$ be its exceptional divisor.
\end{setting}

Casagrande studies $\pi$ when $X$ is a fourfold and $Z_1,Z_2$ are surfaces meeting transversely in one point \cite[Lemma~2.10]{Casagrande}.
We extend her results to \cref{set:main}.

Since $Z_1,\dotsc,Z_s$ meet transversely, they have simple normal crossings in the sense of \cite[Definition~2.1]{MasamuraBlowup}.
Because $X$ is smooth and no $Z_i$ contains another, the centers satisfy \cite[Setting~2.2]{MasamuraBlowup} for any ordering, so $\pi$ is a successive blowup in the sense of that paper.
\Cref{lem:basic} below follows directly from \cite{MasamuraBlowup}, and \cref{lem:local} describes $\pi$ near the centers by its coordinate models.

\begin{lem}\label{lem:basic}
	In \cref{set:main}, the following hold.
	\begin{enumerate}
		\item\label{item:basic-smooth} The variety $\tilde X$ is smooth.
			For every $i$, the divisor $E_i$ is the strict transform of $E_i^{(i)}$, and in particular it is a prime divisor.
		\item\label{item:basic-basis} The following equalities hold:
			\[
				N^1(\tilde X)=\pi^*N^1(X)\oplus\bigoplus_{i=1}^s\RR[E_i],
				\qquad
				K_{\tilde X}-\sum_{i=1}^s(r_i-1)E_i=\pi^*K_X.
			\]
			In particular, $\rho(\tilde X)=\rho(X)+s$.
	\end{enumerate}
\end{lem}

\begin{proof}
	The smoothness of $\tilde X$ and of the centers $Z_i^{(i-1)}$ is \cite[Lemma~2.3]{MasamuraBlowup}.
	For $h>i$, the center $Z_h^{(h-1)}$ maps onto $Z_h\not\subseteq Z_i$, so it is not contained in the strict transform of $E_i^{(i)}$.
	Hence the pullback of $E_i^{(i)}$ coincides with its strict transform.
	For every $i$, we have $N^1(X^{(i)})=\pi_i^*N^1(X^{(i-1)})\oplus\RR[E_i^{(i)}]$, which gives the decomposition of $N^1(\tilde X)$ in \cref{item:basic-basis}.
	The formula for $K_{\tilde X}$ follows from the canonical bundle formula for each $\pi_i$.
\end{proof}

\begin{lem}\label{lem:local}
	In \cref{set:main}, the following hold.
	\begin{enumerate}
		\item\label{item:local-corner} Let $p\in X$ be a closed point.
			Then there is an isomorphism
			\[
				\pi^{-1}(p)\isom\prod_{i\in I(p)}\PP^{r_i-1}
			\]
			under which $\bigrestr{\Oc_{\tilde X}(E_i)}{\pi^{-1}(p)}\isom\mathrm{pr}_i^*\Oc_{\PP^{r_i-1}}(-1)$ for every $i\in I(p)$.
			Moreover, $E_i$ is disjoint from $\pi^{-1}(p)$ for every $i\notin I(p)$.
		\item\label{item:local-ideal} The equality $\Ic_{Z_i}\cdot\Oc_{\tilde X}=\Oc_{\tilde X}(-E_i)$ holds for every $i$.
			In particular, $\pi^{-1}(Z_i)=E_i$.
		\item\label{item:local-fiber} Let $\bar\pi_i\colon\tilde X\to\Bl_{Z_i}X$ be the morphism over $X$ given by \cref{item:local-ideal} and the universal property of the blowup.
			Then the morphism
			\[
				(\bar\pi_1,\dotsc,\bar\pi_s)\colon\tilde X\longrightarrow\Bl_{Z_1}X\times_X\dotsb\times_X\Bl_{Z_s}X
			\]
			is an isomorphism over $X$, and $E_i=\bar\pi_i^*\bar E_i$ for every $i$.
		\item\label{item:local-indep} Up to isomorphism over $X$, the morphism $\pi$ and the divisors $E_1,\dotsc,E_s$ do not depend on the order of $Z_1,\dotsc,Z_s$.
	\end{enumerate}
\end{lem}

\begin{proof}
	We first describe $\pi$ near a closed point $p\in\bigcup_iZ_i$.
	Let $\pi_{Y_p}\colon\tilde Y_p\to\AA^n$ together with $\phi_p\colon U_p\to\AA^n$ be a coordinate model of $\pi$ at $p$ in the sense of \cite[Definition~3.1]{MasamuraBlowup}, with exceptional divisors $F_1,\dotsc,F_s$.
	Set $\tilde U_p\coloneq\pi^{-1}(U_p)$.
	By \cite[Proposition~3.2]{MasamuraBlowup}, there is a commutative diagram
	\[
	\begin{tikzcd}
		\tilde X\arrow[d,"\pi"']
			& \tilde U_p\arrow[l,hook']\arrow[r,"\tilde\phi_p"]\arrow[d]
			& \tilde Y_p\arrow[d,"\pi_{Y_p}"]\\
		X
			& U_p\arrow[l,hook']\arrow[r,"\phi_p"']
			& \AA^n
	\end{tikzcd}
	\]
	in which both squares are cartesian, and $\restr{E_i}{\tilde U_p}=\tilde\phi_p^*F_i$ for every $i$.
	Moreover, $\pi_{Y_p}$ is the successive blowup of $\AA^n$ along the coordinate subspaces $W_i\coloneq V(x_j\mid j\in J_i)$ for $i\in I(p)$, and $F_i=0$ for $i\notin I(p)$.
	We claim that
	\begin{equation}\label{eq:product}
		\tilde Y_p\isom\prod_{i\in I(p)}\Bl_0\AA^{J_i}\times\AA^{J'},
		\qquad J'\coloneq\{1,\dotsc,n\}\setminus\bigcup_{i\in I(p)}J_i,
	\end{equation}
	where $\AA^J$ is the affine space with coordinates $x_j$ for $j\in J$, and that $F_i$ is the pullback of the exceptional divisor of $\Bl_0\AA^{J_i}$.
	Since the sets $J_i$ are pairwise disjoint, $\AA^n=\prod_{i\in I(p)}\AA^{J_i}\times\AA^{J'}$, and $W_i$ is the product of the origin of $\AA^{J_i}$ with the other factors.
	We prove the claim by induction on the number of blowups.
	Let $i\in I(p)$, and assume that the blowups along the $W_h$ with $h\in I(p)$ and $h<i$ replace each factor $\AA^{J_h}$ by $\Bl_0\AA^{J_h}$.
	Then the strict transform of $W_i$ is still the product of the origin of $\AA^{J_i}$ with the other factors, so blowing it up only replaces $\AA^{J_i}$ by $\Bl_0\AA^{J_i}$.
	This also gives the assertion on $F_i$.

	For \cref{item:local-corner}, we may assume $p\in\bigcup_iZ_i$, since otherwise $I(p)=\emptyset$, the fiber $\pi^{-1}(p)$ is a point, and no $E_i$ meets $\pi^{-1}(p)$.
	Then \eqref{eq:product} gives $\pi_{Y_p}^{-1}(0)\isom\prod_{i\in I(p)}\PP^{r_i-1}$, on which $\Oc(F_i)$ restricts to $\mathrm{pr}_i^*\Oc(-1)$ for $i\in I(p)$.
	Since the right square is cartesian and $\phi_p(p)=0$, the morphism $\tilde\phi_p$ induces an isomorphism $\pi^{-1}(p)\isomarrow\pi_{Y_p}^{-1}(0)$.

	For \cref{item:local-ideal}, both sides are trivial over $X\setminus Z_i$ by \cref{item:local-corner}, so it suffices to work over $U_p$ for $p\in Z_i$.
	There $\Ic_{Z_i}$ is generated by the pullbacks of the $x_j$ with $j\in J_i$, and $(x_j\mid j\in J_i)\cdot\Oc_{\tilde Y_p}=\Oc_{\tilde Y_p}(-F_i)$ by \eqref{eq:product}.
	Pulling back along $\tilde\phi_p$ gives $\Ic_{Z_i}\cdot\Oc_{\tilde U_p}=\Oc_{\tilde U_p}(-E_i)$.

	We next prove \cref{item:local-fiber}.
	Set $\bar X\coloneq\Bl_{Z_1}X\times_X\dotsb\times_X\Bl_{Z_s}X$ and $\bar\pi\coloneq(\bar\pi_1,\dotsc,\bar\pi_s)\colon\tilde X\to\bar X$.
	Since $\bar\pi$ is the identity of $X$ over $X\setminus\bigcup_iZ_i$, it suffices to show that $\bar\pi$ is an isomorphism over $U_p$ for every closed point $p\in\bigcup_iZ_i$.
	Set
	\[
		\bar U_p\coloneq\bar X\times_XU_p,
		\qquad
		\bar Y_p\coloneq\prod_{i\in I(p)}{}_{\AA^n}\,\Bl_{W_i}\AA^n,
	\]
	where the product is the fiber product over $\AA^n$.
	Recall that $Z_i\cap U_p=\phi_p^{-1}(W_i)$ for $i\in I(p)$ and $Z_i\cap U_p=\emptyset$ for $i\notin I(p)$.
	Since $\phi_p$ is \'etale and blowups commute with flat base change, we obtain a commutative diagram
	\[
	\begin{tikzcd}
		\bar X\arrow[d]
			& \bar U_p\arrow[l,hook']\arrow[r]\arrow[d]
			& \bar Y_p\arrow[d]\\
		X
			& U_p\arrow[l,hook']\arrow[r,"\phi_p"']
			& \AA^n
	\end{tikzcd}
	\]
	in which both squares are cartesian.
	Let $\bar\pi_{Y_p}\colon\tilde Y_p\to\bar Y_p$ be the morphism over $\AA^n$ given by the universal property of the blowups.
	Since the sets $J_i$ are disjoint, $\bar Y_p$ is the right-hand side of \eqref{eq:product}, and $\bar\pi_{Y_p}$ is an isomorphism.
	By the uniqueness in the universal property, $\bar\pi$ over $U_p$ is the base change of $\bar\pi_{Y_p}$ along $\phi_p$.
	Hence $\bar\pi$ is an isomorphism over $U_p$.
	Finally, $E_i=\bar\pi_i^*\bar E_i$, since both are effective Cartier divisors with ideal $\Ic_{Z_i}\cdot\Oc_{\tilde X}$ by \cref{item:local-ideal}.

	We obtain \cref{item:local-indep} from \cref{item:local-fiber}, applied to every ordering of $Z_1,\dotsc,Z_s$, since neither the fiber product nor the divisors $\bar\pi_i^*\bar E_i$ depend on the order.
\end{proof}

\begin{lem}\label{lem:relample}
	In \cref{set:main}, the divisor $-K_{\tilde X}$ is $\pi$-ample.
\end{lem}

\begin{proof}
	It suffices to show that $-K_{\tilde X}$ is ample on $\pi^{-1}(p)$ for every closed point $p\in X$.
	By \cref{lem:basic}~\cref{item:basic-basis} and \cref{lem:local}~\cref{item:local-corner}, the restriction of $\Oc_{\tilde X}(-K_{\tilde X})$ to $\pi^{-1}(p)\isom\prod_{i\in I(p)}\PP^{r_i-1}$ is isomorphic to $\bigotimes_{i\in I(p)}\mathrm{pr}_i^*\Oc_{\PP^{r_i-1}}(r_i-1)$.
	This is ample, since $r_i\ge2$ for every $i$.
\end{proof}

By transversality, a curve on $X$ lifts to $\tilde X$ by choosing a normal direction along it independently for each center containing it, and this fixes its intersection numbers with $E_i$.

\begin{lem}\label{lem:lift}
	In \cref{set:main}, let $\bar C$ be a smooth projective curve, and let $g\colon\bar C\to X$ be a morphism whose image is a curve $C$.
	For $i$ with $C\subseteq Z_i$, let $g^*N_{Z_i/X}^\vee\to\Lc_i$ be an invertible quotient.
	Then there is a morphism $\tilde g\colon\bar C\to\tilde X$ with $\pi\circ\tilde g=g$ such that, for every $i$,
	\[
		\deg\tilde g^*E_i
		\begin{cases}
			=-\deg\Lc_i&\text{if }C\subseteq Z_i,\\
			\ge1&\text{if }C\not\subseteq Z_i\text{ and }C\cap Z_i\ne\emptyset,\\
			=0&\text{if }C\cap Z_i=\emptyset.
		\end{cases}
	\]
\end{lem}

\begin{proof}
	By \cref{lem:local}~\cref{item:local-fiber}, it suffices to construct, for every $i$, a morphism $g_i\colon\bar C\to\Bl_{Z_i}X$ with $\pi'_i\circ g_i=g$ such that $\deg g_i^*\bar E_i$ has the required property.
	Indeed, the morphism $\tilde g$ with $\bar\pi_i\circ\tilde g=g_i$ for every $i$ then satisfies $\tilde g^*E_i=g_i^*\bar E_i$, since $E_i=\bar\pi_i^*\bar E_i$.

	Assume first that $C\not\subseteq Z_i$.
	Then the ideal $g^{-1}\Ic_{Z_i}\cdot\Oc_{\bar C}$ is nonzero, and hence invertible, since $\bar C$ is a smooth curve.
	By the universal property of the blowup, there is a morphism $g_i\colon\bar C\to\Bl_{Z_i}X$ with $\pi'_i\circ g_i=g$ and
	\[
		g_i^*\Oc(-\bar E_i)=g^{-1}\Ic_{Z_i}\cdot\Oc_{\bar C}.
	\]
	Hence $\deg g_i^*\bar E_i$ is the degree of the effective divisor on $\bar C$ defined by $g^{-1}\Ic_{Z_i}\cdot\Oc_{\bar C}$.
	Its support is $g^{-1}(Z_i)$, so this degree is at least $1$ if $C\cap Z_i\ne\emptyset$ and $0$ otherwise.

	Assume next that $C\subseteq Z_i$.
	Then $g$ factors through $Z_i$.
	We have $\bar E_i\isom\PP_{Z_i}(N_{Z_i/X}^\vee)$ and $\Oc_{\bar E_i}(\bar E_i)\isom\Oc(-1)$.
	Hence the invertible quotient $g^*N_{Z_i/X}^\vee\to\Lc_i$ defines a morphism $g_i\colon\bar C\to\bar E_i\subseteq\Bl_{Z_i}X$ over $g$ with $g_i^*\Oc_{\bar E_i}(\bar E_i)\isom\Lc_i^{-1}$.
	Thus $\deg g_i^*\bar E_i=-\deg\Lc_i$.
\end{proof}

\section{Fano blowups of Grassmannians}\label{sec:grass}

In this section, we consider the successive blowup of a Grassmannian $\Gr(d,V)$ along sub-Grassmannians $\Gr(d,W_i)$ that meet transversely.
We characterize, in terms of the subspaces $W_i$, when the blowup is Fano, and we apply the result to the Grassmannians of dimension four.

\subsection{Setting}\label{ss:subgrass}

In this subsection, we collect basic facts on sub-Grassmannians and fix the setting of this section.

We use the notation for Grassmannians of \cref{ss:conventions}.
For a subspace $W\subseteq V$, we write $W^\perp\subseteq V^\vee$ for the space of linear forms vanishing on $W$.

\begin{lem}\label{lem:subgrass}
	Let $V$ be a vector space of dimension $n+1$, let $1\le d\le n$, and set $X\coloneq\Gr(d,V)$.
	Let $W_1,\dotsc,W_s\subsetneq V$ be subspaces of dimension at least $d$, let $c_i$ be the codimension of $W_i$ in $V$, and set $Z_i\coloneq\Gr(d,W_i)$.
	Then the following hold.
	\begin{enumerate}
		\item\label{item:subgrass-zero} The subvariety $Z_i$ is the zero scheme of the section $\sigma_i\in H^0(X,\Sc^\vee\otimes(V/W_i))$ given by the composite
			\[
				\Sc\longrightarrow V\otimes\Oc_X\longrightarrow(V/W_i)\otimes\Oc_X.
			\]
			It is smooth of codimension $dc_i$, and $N_{Z_i/X}\isom\restr{\Sc^\vee}{Z_i}\otimes(V/W_i)$.
		\item\label{item:subgrass-transverse} The subvarieties $Z_1,\dotsc,Z_s$ meet transversely if and only if, for every subset $I\subseteq\{1,\dotsc,s\}$ such that $\bigcap_{i\in I}W_i$ has dimension at least $d$, the codimension of $\bigcap_{i\in I}W_i$ in $V$ is $\sum_{i\in I}c_i$.
	\end{enumerate}
\end{lem}

\begin{proof}
	Let $g\colon T\to X$ be a morphism from a $k$-scheme $T$, and set $\Sc_T\coloneq g^*\Sc\subseteq V\otimes\Oc_T$.
	By the universal property of Grassmannians, $g$ factors through $Z_i=\Gr(d,W_i)$ if and only if $\Sc_T\subseteq W_i\otimes\Oc_T$ and $(W_i\otimes\Oc_T)/\Sc_T$ is locally free.
	The second condition follows from the first, because $(W_i\otimes\Oc_T)/\Sc_T$ is then the kernel of the surjection $(V\otimes\Oc_T)/\Sc_T\to(V/W_i)\otimes\Oc_T$ of locally free sheaves.
	Since $g^*\sigma_i\in H^0(T,\Sc_T^\vee\otimes(V/W_i))$ is given by the composite $\Sc_T\to V\otimes\Oc_T\to(V/W_i)\otimes\Oc_T$, the first condition means that $g^*\sigma_i=0$, that is, $g$ factors through the zero scheme of $\sigma_i$.
	Therefore $Z_i$ and the zero scheme of $\sigma_i$ coincide as closed subschemes of $X$.
	Since $Z_i=\Gr(d,W_i)$ is smooth of dimension $d(n+1-c_i-d)$, its codimension in $X$ is $dc_i$, which is the rank of $\Sc^\vee\otimes(V/W_i)$.
	Thus $\sigma_i$ is a regular section, and $N_{Z_i/X}\isom\restr{\Sc^\vee}{Z_i}\otimes(V/W_i)$.
	This proves \cref{item:subgrass-zero}.

	For \cref{item:subgrass-transverse}, first assume that $Z_1,\dotsc,Z_s$ meet transversely.
	Let $I\subseteq\{1,\dotsc,s\}$ be a subset, and let $P\subseteq\bigcap_{i\in I}W_i$ be a $d$-dimensional subspace.
	Then $[P]\in Z_i$ for every $i\in I$, so by transversality at $[P]$, the intersection $\bigcap_{i\in I}Z_i=\Gr(d,\bigcap_{i\in I}W_i)$ has codimension $\sum_{i\in I}dc_i$ near $[P]$.
	Thus $\bigcap_{i\in I}W_i$ has codimension $\sum_{i\in I}c_i$ in $V$.

	Conversely, assume the condition in \cref{item:subgrass-transverse}.
	Let $p=[P]\in\bigcup_iZ_i$ be a closed point.
	We construct coordinates centered at $p$ as in \cref{defi:transverse}.
	Choose a complement $V_0$ of $P$ in $V$.
	Sending a linear map $f\colon P\to V_0$ to its graph $\sset{v+f(v)}{v\in P}$ gives an open immersion $\Hom(P,V_0)\to X$ onto the neighborhood $U$ of $p$ consisting of the $d$-dimensional subspaces $P'$ with $P'\cap V_0=0$, and $p$ corresponds to $f=0$.
	For $i\in I(p)$, we have $P\subseteq W_i$, so the graph of $f$ is contained in $W_i$ if and only if $f(P)\subseteq V_0\cap W_i$.
	Thus $Z_i\cap U$ is the linear subspace $\Hom(P,V_0\cap W_i)$ of $\Hom(P,V_0)=P^\vee\otimes V_0$.
	Let $L_i\subseteq V_0^\vee$ be the space of linear forms on $V_0$ vanishing on $V_0\cap W_i$.
	Then $Z_i\cap U$ is defined by the linear forms in $P\otimes L_i$.
	Hence it suffices to show that the sum $\sum_{i\in I(p)}L_i\subseteq V_0^\vee$ is direct, because then so is the sum $\sum_{i\in I(p)}(P\otimes L_i)\subseteq P\otimes V_0^\vee$, and extending bases of these subspaces to a basis of $P\otimes V_0^\vee$ gives the desired coordinates on $U$.
	Since $P\subseteq W_i$ for $i\in I(p)$, we have $W_i=P\oplus(V_0\cap W_i)$ and $\bigcap_{i\in I(p)}W_i=P\oplus\bigcap_{i\in I(p)}(V_0\cap W_i)$.
	Hence $\dim L_i=c_i$, and by the condition, $\bigcap_{i\in I(p)}(V_0\cap W_i)$ has codimension $\sum_{i\in I(p)}c_i$ in $V_0$.
	This codimension is the dimension of $\sum_{i\in I(p)}L_i$, so the sum $\sum_{i\in I(p)}L_i$ is direct.
	This proves \cref{item:subgrass-transverse}.
\end{proof}

We now fix the setting of this section.

\begin{setting}\label{set:grass}
	In \cref{set:main}, assume that $X=\Gr(d,V)$ for a vector space $V$ of dimension $n+1$ and an integer $1\le d\le n$, and that $Z_i=\Gr(d,W_i)$ for subspaces $W_1,\dotsc,W_s\subsetneq V$ of dimension at least $d$.
	Let $c_i$ be the codimension of $W_i$ in $V$, so that $r_i=dc_i$.
	Set $H\coloneq\pi^*\Oc_X(1)$.
	For a subset $I\subseteq\{1,\dotsc,s\}$, set $W_I\coloneq\bigcap_{i\in I}W_i$, where $W_\emptyset\coloneq V$.
\end{setting}

Since $-K_X=\Oc_X(n+1)$, \cref{lem:basic}~\cref{item:basic-basis} gives
\begin{equation}\label{eq:grass-antican}
	-K_{\tilde X}=(n+1)H-\sum_{i=1}^s(r_i-1)E_i.
\end{equation}

Whether $\tilde X$ is Fano depends on the following number (\cref{thm:grass}).

\begin{defi}\label{defi:theta}
	In \cref{set:grass}, set
	\[
		\theta\coloneq\max\biggl\{\frac1{\max\{1,d-\dim W_I\}}\sum_{i\in I}(r_i-1)\biggm|I\subseteq\{1,\dotsc,s\}\biggr\}.
	\]
\end{defi}

Intuitively, the term for $I$ compares the sum of the coefficients $r_i-1$ of $E_i$ with $i\in I$ in \eqref{eq:grass-antican} to the degree $\max\{1,d-\dim W_I\}$ of curves on $X$ that meet $Z_i$ for every $i\in I$.

\subsection{Positivity of the anticanonical divisor}\label{ss:grass-positive}

In this subsection, we show that $-K_{\tilde X}-(n+1-\theta)H$ is nef, so that $\tilde X$ is Fano if $\theta<n+1$.

We regard a linear form $f\in V^\vee$ as a section of $\Sc^\vee$ via the natural surjection $V^\vee\otimes\Oc_X\to\Sc^\vee$.
For every $i$, we write $s_{E_i}$ for the canonical section of $\Oc_{\tilde X}(E_i)$.

\begin{lem}\label{lem:grass-sections}
	In \cref{set:grass}, for every $i$, there is a morphism $\varphi_i\colon W_i^\perp\otimes\Oc_{\tilde X}\to\pi^*\Sc^\vee(-E_i)$ which is nonzero at every point of $\tilde X$ and makes the diagram
	\[
	\begin{tikzcd}
		W_i^\perp\otimes\Oc_{\tilde X}\arrow[r,hook]\arrow[d,"\varphi_i"']
			& V^\vee\otimes\Oc_{\tilde X}\arrow[d]\\
		\pi^*\Sc^\vee(-E_i)\arrow[r,hook,"s_{E_i}"']
			& \pi^*\Sc^\vee
	\end{tikzcd}
	\]
	commute, where the right vertical map is the pullback of $V^\vee\otimes\Oc_X\to\Sc^\vee$.
\end{lem}

\begin{proof}
	Via $W_i^\perp=(V/W_i)^\vee$, the composite of $W_i^\perp\otimes\Sc\hookrightarrow V^\vee\otimes\Sc$ and the natural map $V^\vee\otimes\Sc\to\Oc_X$ is the dual of $\sigma_i\colon\Oc_X\to\Sc^\vee\otimes(V/W_i)$ in \cref{lem:subgrass}~\cref{item:subgrass-zero}, so its image is the ideal $\Ic_{Z_i}$ of the zero scheme $Z_i$ of $\sigma_i$.
	Hence the image of its pullback $W_i^\perp\otimes\pi^*\Sc\to\Oc_{\tilde X}$ is $\Ic_{Z_i}\cdot\Oc_{\tilde X}=\Oc_{\tilde X}(-E_i)$ by \cref{lem:local}~\cref{item:local-ideal}, and we obtain a commutative diagram
	\[
	\begin{tikzcd}
		W_i^\perp\otimes\pi^*\Sc\arrow[r,hook]\arrow[d,two heads]
			& V^\vee\otimes\pi^*\Sc\arrow[d]\\
		\Oc_{\tilde X}(-E_i)\arrow[r,hook,"s_{E_i}"']
			& \Oc_{\tilde X}
	\end{tikzcd}
	\]
	By the adjunction $\Hom(\Ec\otimes\pi^*\Sc,\Fc)\isom\Hom(\Ec,\pi^*\Sc^\vee\otimes\Fc)$, this diagram gives the diagram in the statement.
	Since the left vertical map in the above diagram is surjective, $\varphi_i$ is nonzero at every point.
\end{proof}

For $f\in W_i^\perp$, we simply write $\varphi_i(f)$ for $\varphi_i(f\otimes1)\in H^0(\tilde X,\pi^*\Sc^\vee(-E_i))$.
By \cref{lem:grass-sections}, it satisfies $s_{E_i}\varphi_i(f)=\pi^*f$.
We now prove the positivity of $-K_{\tilde X}$ by comparing these sections with the Harder--Narasimhan filtration of $\pi^*\Sc^\vee$ on curves.

\begin{prop}\label{prop:grass-positive}
	In \cref{set:grass}, the $\QQ$-divisor
	\[
		-K_{\tilde X}-(n+1-\theta)H=\theta H-\sum_{i=1}^s(r_i-1)E_i
	\]
	is nef.
	In particular, $\tilde X$ is Fano if $\theta<n+1$.
\end{prop}

\begin{proof}
	The equality follows from \eqref{eq:grass-antican}.
	We first deduce the last assertion from the nefness of this $\QQ$-divisor.
	Assume that $\theta<n+1$.
	By \cref{lem:relample}, the divisor $-K_{\tilde X}$ is $\pi$-ample, so $-K_{\tilde X}+tH$ is ample for $t\gg0$.
	Hence
	\[
		-K_{\tilde X}=\frac{n+1-\theta}{t+n+1-\theta}(-K_{\tilde X}+tH)+\frac{t}{t+n+1-\theta}(-K_{\tilde X}-(n+1-\theta)H)
	\]
	is ample.

	It remains to show that $\theta H-\sum_i(r_i-1)E_i$ is nef.
	Let $C\subseteq\tilde X$ be a curve, let $\nu\colon\bar C\to\tilde X$ be its normalization, and set $\Sc_{\bar C}\coloneq\nu^*\pi^*\Sc$.
	Since $\deg\Sc^\vee_{\bar C}=H\cdot C$, it suffices to show that
	\[
		\sum_{i=1}^s(r_i-1)E_i\cdot C\le\theta\deg\Sc^\vee_{\bar C}.
	\]

	For a nonzero torsion-free sheaf $G$ on $\bar C$, we write $\mu(G)\coloneq\deg G/\rank G$ for its slope and $G(-E_i)\coloneq G\otimes\nu^*\Oc_{\tilde X}(-E_i)$.
	Let
	\[
		0=F_0\subsetneq F_1\subsetneq\dotsb\subsetneq F_t=\Sc^\vee_{\bar C}
	\]
	be the Harder--Narasimhan filtration of $\Sc^\vee_{\bar C}$.
	Each quotient $F_j/F_{j-1}$ is semistable, and $\mu_j\coloneq\mu(F_j/F_{j-1})$ satisfies $\mu_1>\dotsb>\mu_t$.
	For every $i$, since $\varphi_i$ in \cref{lem:grass-sections} is nonzero at every point, its pullback
	\[
		\nu^*\varphi_i\colon W_i^\perp\otimes\Oc_{\bar C}\longrightarrow\Sc^\vee_{\bar C}(-E_i)
	\]
	is nonzero.
	Let $j(i)$ be the smallest integer $j$ such that the image of $\nu^*\varphi_i$ is contained in $F_j(-E_i)$.

	\begin{claim}\label{claim:grass-HN}
		For every $i$, we have $E_i\cdot C\le\mu_{j(i)}$.
	\end{claim}

	\begin{proof}[Proof of \cref{claim:grass-HN}]
		By the choice of $j(i)$, the composite
		\[
			W_i^\perp\otimes\Oc_{\bar C}\xrightarrow{\nu^*\varphi_i}F_{j(i)}(-E_i)\longrightarrow(F_{j(i)}/F_{j(i)-1})(-E_i)
		\]
		is nonzero.
		Let $G$ be its image.
		Since $G$ is a nonzero subsheaf of $(F_{j(i)}/F_{j(i)-1})(-E_i)$, which is semistable of slope $\mu_{j(i)}-E_i\cdot C$, it is torsion-free and satisfies $\mu(G)\le\mu_{j(i)}-E_i\cdot C$.
		Since $G$ is a quotient of $W_i^\perp\otimes\Oc_{\bar C}$, we have $\mu(G)\ge0$.
		Hence $E_i\cdot C\le\mu_{j(i)}$.
	\end{proof}

	\begin{claim}\label{claim:grass-partial}
		For $1\le j\le t$, we have
		\[
			a_j\coloneq\sum_{j(i)\le j}(r_i-1)\le\theta\rank F_j.
		\]
	\end{claim}

	\begin{proof}[Proof of \cref{claim:grass-partial}]
		Let $L_j\subseteq V^\vee$ be the kernel of the composite $V^\vee\to H^0(\bar C,\Sc^\vee_{\bar C})\to H^0(\bar C,\Sc^\vee_{\bar C}/F_j)$.
		Since $\Sc^\vee_{\bar C}/F_j$ is a quotient of $V^\vee\otimes\Oc_{\bar C}$, the image of $V^\vee$ generates $\Sc^\vee_{\bar C}/F_j$, which has rank $d-\rank F_j$, so $L_j$ has codimension at least $d-\rank F_j$.
		Hence $L_j^\perp\subseteq V$ has dimension at least $d-\rank F_j$.
		Set $I_j\coloneq\sset{i}{j(i)\le j}$, and let $i\in I_j$.
		For every $f\in W_i^\perp$, the section $\nu^*\varphi_i(f)$ lies in $H^0(\bar C,F_j(-E_i))$, so $\nu^*\pi^*f=\nu^*s_{E_i}\cdot\nu^*\varphi_i(f)$ lies in $H^0(\bar C,F_j)$.
		Thus $W_i^\perp\subseteq L_j$, that is, $L_j^\perp\subseteq W_i$.
		Therefore $L_j^\perp\subseteq W_{I_j}$, so $d-\dim W_{I_j}\le\rank F_j$.
		Since $\rank F_j\ge1$, \cref{defi:theta} gives $a_j\le\theta\max\{1,d-\dim W_{I_j}\}\le\theta\rank F_j$.
	\end{proof}

	Since $\Sc^\vee_{\bar C}$ is a quotient of $V^\vee\otimes\Oc_{\bar C}$, so is $\Sc^\vee_{\bar C}/F_{t-1}$, and hence $\mu_t\ge0$.
	Set $a_0=\mu_{t+1}\coloneq0$.
	Since $r_i-1\ge0$ and $\mu_j-\mu_{j+1}\ge0$ for $1\le j\le t$, \cref{claim:grass-HN,claim:grass-partial} give
	\begin{align*}
		\sum_{i=1}^s(r_i-1)E_i\cdot C
		&\le\sum_{j=1}^t(a_j-a_{j-1})\mu_j
		=\sum_{j=1}^ta_j(\mu_j-\mu_{j+1})\\
		&\le\theta\sum_{j=1}^t\rank F_j\,(\mu_j-\mu_{j+1})
		=\theta\sum_{j=1}^t\rank(F_j/F_{j-1})\,\mu_j
		=\theta\deg\Sc^\vee_{\bar C}. \qedhere
	\end{align*}
\end{proof}

\subsection{Rational curves meeting the centers}\label{ss:grass-curves}

In this subsection, we construct rational curves on $X$ that meet many of the centers $Z_i$, and we bound the anticanonical degree of their lifts to $\tilde X$.

We first apply \cref{lem:lift} to rational curves on $X$ on which $\Sc^\vee$ splits in a simple way.

\begin{lem}\label{lem:grass-lift}
	In \cref{set:grass}, let $g\colon\PP^1\to X$ be a morphism such that
	\[
		g^*\Sc^\vee\isom\Oc_{\PP^1}(1)^{\oplus a}\oplus\Oc_{\PP^1}^{\oplus(d-a)}
	\]
	for an integer $1\le a\le d$.
	Let $I$ be a set of indices $i$ with $g(\PP^1)\cap Z_i\ne\emptyset$.
	Then there is a morphism $\tilde g\colon\PP^1\to\tilde X$ with $\pi\circ\tilde g=g$ and
	\[
		\deg\tilde g^*(-K_{\tilde X})\le(n+1)a-\sum_{i\in I}(r_i-1).
	\]
	In particular, if $\sum_{i\in I}(r_i-1)\ge(n+1)a$, then $\tilde X$ is not Fano.
\end{lem}

\begin{proof}
	The image $C$ of $g$ is a curve, because $g^*\Oc_X(1)\isom\det g^*\Sc^\vee$ has degree $a>0$.
	If $C\subseteq Z_i$, then \cref{lem:subgrass}~\cref{item:subgrass-zero} gives
	\[
		g^*N_{Z_i/X}^\vee\isom(\Oc_{\PP^1}(-1)^{\oplus a}\oplus\Oc_{\PP^1}^{\oplus(d-a)})^{\oplus c_i},
	\]
	and the projection to one summand gives a quotient $\Lc_i\coloneq\Oc_{\PP^1}(-1)$.
	With these quotients, \cref{lem:lift} gives a morphism $\tilde g\colon\PP^1\to\tilde X$ with $\pi\circ\tilde g=g$ such that, for every $i$,
	\[
		\deg\tilde g^*E_i
		\begin{cases}
			=-\deg\Lc_i=1&\text{if }C\subseteq Z_i,\\
			\ge1&\text{if }C\not\subseteq Z_i\text{ and }C\cap Z_i\ne\emptyset,\\
			=0&\text{if }C\cap Z_i=\emptyset.
		\end{cases}
	\]
	In particular, $\deg\tilde g^*E_i\ge1$ for every $i\in I$, and $\deg\tilde g^*E_i\ge0$ for every $i$.
	Since $\deg\tilde g^*H=a$, the inequality follows from \eqref{eq:grass-antican}.
	Since $\deg\tilde g^*(-K_{\tilde X})$ is a positive multiple of $-K_{\tilde X}\cdot\tilde g(\PP^1)$, this proves the last assertion.
\end{proof}

By a line on $X$ we mean a curve of degree one with respect to $\Oc_X(1)$.
The restriction of $\Sc^\vee$ to a line is $\Oc_{\PP^1}(1)\oplus\Oc_{\PP^1}^{\oplus(d-1)}$, so \cref{lem:grass-lift} applies to the inclusion of a line with $a=1$.

\begin{lem}\label{lem:grass-line}
	In \cref{set:grass}, let $I$ be a set of indices with
	\[
		\dim W_I\ge d-1,\qquad\sum_{i\in I}(c_i-1)\le2(n-d).
	\]
	Then there is a line in $X$ that meets $Z_i$ for every $i\in I$.
\end{lem}

\begin{proof}
	Fix a subspace $V_0\subseteq W_I$ of dimension $d-1$ and a generator $v_0$ of $\bigwedge^{d-1}V_0$, and put $m\coloneq n+1-d$.
	As every $v_0\wedge v$ is decomposable, the map $V/V_0\hookrightarrow\bigwedge^dV$, $v\mapsto v_0\wedge v$, embeds $\PP^m\coloneq\PP((V/V_0)^\vee)$ linearly into $X\subseteq\PP((\bigwedge^dV)^\vee)$.
	For $i\in I$, the linear subspace $\Pi_i\coloneq\PP((W_i/V_0)^\vee)$ of $\PP^m$ has codimension $c_i$.
	The same argument with $W_i$ in place of $V$ shows that this embedding maps $\Pi_i$ into $Z_i=\Gr(d,W_i)$.
	Hence it suffices to find a line in $\PP^m$ that meets $\Pi_i$ for every $i\in I$.

	The lines meeting $\Pi_i$ form the special Schubert variety $\Sigma_i\subseteq\Gr(2,V/V_0)$, whose class in $A^{c_i-1}(\Gr(2,V/V_0))$ is $\sigma_{c_i-1}$.
	Hence it suffices to show that $\prod_{i\in I}\sigma_{c_i-1}\ne0$ in the Chow ring $A^*(\Gr(2,V/V_0))$, since then $\bigcap_{i\in I}\Sigma_i\ne\emptyset$.

	The Schubert classes $\sigma_\lambda$, indexed by sequences of integers $\lambda=(\lambda_1,\lambda_2)$ with $m-1\ge\lambda_1\ge\lambda_2\ge0$, form a $\ZZ$-basis of $A^*(\Gr(2,V/V_0))$ \cite[Corollary~4.7]{EisenbudHarris}, and by the Pieri formula \cite[Proposition~4.9]{EisenbudHarris} every product $\sigma_\lambda\cdot\sigma_b$ with an integer $0\le b\le m-1$ is a nonnegative combination of them.
	For $0\le a\le2m-2$, set $\lambda(a)\coloneq(a,0)$ if $a\le m-1$ and $\lambda(a)\coloneq(m-1,a-m+1)$ if $a\ge m-1$.
	Set $a_I\coloneq\sum_{i\in I}(c_i-1)$, so that $a_I\le2(n-d)=2m-2$ by assumption.
	We show by induction on $\card I$ that $\sigma_{\lambda(a_I)}$ occurs with positive coefficient in $\prod_{i\in I}\sigma_{c_i-1}$, which implies $\prod_{i\in I}\sigma_{c_i-1}\ne0$.
	If $I=\emptyset$, then $\sigma_{\lambda(0)}=1$ is the empty product.
	Otherwise, choose $i\in I$ and set $I'\coloneq I\setminus\{i\}$.
	By the induction hypothesis and the nonnegativity above, it suffices to show that $\sigma_{\lambda(a_I)}$ occurs with positive coefficient in $\sigma_{\lambda(a_{I'})}\cdot\sigma_{c_i-1}$.
	Let $\lambda(a_I)\eqcolon(\lambda_1,\lambda_2)$ and $\lambda(a_{I'})\eqcolon(\lambda'_1,\lambda'_2)$.
	Then $\lambda_1+\lambda_2=\lambda'_1+\lambda'_2+(c_i-1)$, and $\lambda'_1\le\lambda_1$ and $\lambda'_2\le\lambda_2$ since $a_{I'}\le a_I$.
	Moreover, $\lambda_2\le\lambda'_1$ since $c_i-1\le m-1$ and $a_I\le2m-2$, so the claim follows from the Pieri formula.
\end{proof}

We next construct rational curves as in \cref{lem:grass-lift} with $a=2$ for $d=2$.

\begin{lem}\label{lem:grass-conic}
	In \cref{set:grass}, assume that $d=2$.
	Let $I$ be a set of indices, and assume that one of the following holds.
	\begin{enumerate}
		\item\label{item:conic-quadric} $r_i=2$ for every $i\in I$, and $\card I\le3n-1$.
			Moreover, $\card{I_0}<\card I/2$ for every subset $I_0\subseteq I$ with $W_{I_0}\ne0$.
		\item\label{item:conic-graph} $\card I=3$, $r_i=n+1$ for every $i\in I$, and $W_i\cap W_j=0$ for distinct $i,j\in I$.
	\end{enumerate}
	Then there is a morphism $g\colon\PP^1\to X$ with $g^*\Sc^\vee\isom\Oc_{\PP^1}(1)^{\oplus2}$ whose image meets $Z_i$ for every $i\in I$.
\end{lem}

\begin{proof}
	We have isomorphisms
	\begin{align*}
		\Hom(\Oc_{\PP^1}(-1)^{\oplus2},V\otimes\Oc_{\PP^1})
		&\isom\Hom_k(k^{\oplus2},V)\otimes H^0(\PP^1,\Oc_{\PP^1}(1))\\
		&\isom\Hom_k(k^{\oplus2}\otimes H^0(\PP^1,\Oc_{\PP^1}(1))^\vee,V)\\
		&\isom\Hom_k(k^{\oplus2}\otimes k^{\oplus2},V),
	\end{align*}
	where we take $\PP^1=\PP((k^{\oplus2})^\vee)$ for the last one.
	Hence a morphism $\varphi\colon\Oc_{\PP^1}(-1)^{\oplus2}\to V\otimes\Oc_{\PP^1}$ corresponds to the linear map $y\otimes x\mapsto{}^tyMx$ for a unique $2\times2$ matrix $M$ with entries in $V$, where $y,x\in k^{\oplus2}$ are column vectors.
	For $x\ne0$, the restriction of this map to $k^{\oplus2}\otimes kx$ is $\varphi\otimes k([x])$ via the tautological inclusion $\Oc_{\PP^1}(-1)\hookrightarrow k^{\oplus2}\otimes\Oc_{\PP^1}$.
	Therefore, if ${}^tyMx\ne0$ for all $y,x\in k^{\oplus2}\setminus\{0\}$, then $\varphi\otimes k([x])$ is injective for every $x\ne0$, so $\varphi$ defines a morphism $g\colon\PP^1\to X$ with $g^*\Sc\isom\Oc_{\PP^1}(-1)^{\oplus2}$.
	Hence it suffices to find such $M$ that, for every $i\in I$, the entries of $Mx$ lie in $W_i$ for some $x\ne0$, since the image of $g$ then meets $Z_i$.

	Assume first that \cref{item:conic-quadric} holds.
	For $i\in I$, since $r_i=2$, the subspace $W_i$ is a hyperplane in $V$, so $W_i^\perp\subseteq V^\vee$ is one-dimensional.
	Choose a subset $J\subseteq I$ with $\card J\le n-1$ and $\card{I\setminus J}\le2n$.
	Since $W_J$ has dimension at least two, it contains linearly independent vectors $v_1,v_2$.
	Then every element $q$ of the $(2n+1)$-dimensional subspace $v_1V+v_2V\subseteq\Sym^2V$, regarded as a homogeneous quadratic polynomial on $V^\vee$, satisfies $q(W_i^\perp)=0$ for $i\in J$.
	Since the conditions $q(W_i^\perp)=0$ for $i\in I\setminus J$ are at most $2n$ linear conditions, there is a nonzero element $q=v_1u_1+v_2u_2$ of this subspace with $u_1,u_2\in V$ such that $q(W_i^\perp)=0$ for every $i\in I$.
	Then $q=\det M$ for
	\[
		M\coloneq\begin{pmatrix}v_1&-v_2\\u_2&u_1\end{pmatrix}.
	\]

	For $f\in V^\vee$, let $f(M)$ be the matrix obtained by applying $f$ to the entries of $M$, so that $\det f(M)=q(f)$.
	For $i\in I$ and $f\in W_i^\perp\setminus\{0\}$, we have $\det f(M)=0$, so there is $x\ne0$ with $f(M)x=0$, that is, the entries of $Mx$ lie in $W_i=\ker f$.
	We show that ${}^tyMx\ne0$ for all $y,x\in k^{\oplus2}\setminus\{0\}$.
	Suppose, to the contrary, that ${}^tyMx=0$ for some $y,x\in k^{\oplus2}\setminus\{0\}$.
	Choose $S,T\in\GL_2(k)$ such that the first row of $S$ is ${}^ty$ and the first column of $T$ is $x$.
	Then the upper left entry of $SMT$ is $0$, so $\det(ST)q=\det(SMT)=-vv'$ for the off-diagonal entries $v,v'\in V$ of $SMT$.
	Since $q\ne0$, the vectors $v$ and $v'$ are nonzero.
	Since $q(W_i^\perp)=0$ for $i\in I$, the sets $I_0\coloneq\sset{i\in I}{v\in W_i}$ and $I_0'\coloneq\sset{i\in I}{v'\in W_i}$ cover $I$.
	Since $v\in W_{I_0}$ and $v'\in W_{I_0'}$, this contradicts the assumption.
	This proves the assertion in the case \cref{item:conic-quadric}.

	Finally, assume that \cref{item:conic-graph} holds.
	By renumbering, we may assume that $I=\{1,2,3\}$.
	Each $W_i$ has dimension $(n+1)/2$ since $r_i=n+1$, so $W_1\cap W_2=0$ gives $V=W_1\oplus W_2$.
	Since $W_3\cap W_2=0$ and $\dim W_3=\dim W_1$, the subspace $W_3$ is the graph $\sset{w+\alpha(w)}{w\in W_1}$ of a linear map $\alpha\colon W_1\to W_2$, which is injective since $W_3\cap W_1=0$.
	Choose linearly independent vectors $w_1,w_2\in W_1$, and set
	\[
		M\coloneq\begin{pmatrix}w_1&\alpha(w_1)\\w_2&\alpha(w_2)\end{pmatrix}.
	\]
	For $y,x\in k^{\oplus2}\setminus\{0\}$, the vector $w\coloneq y_1w_1+y_2w_2\in W_1$ is nonzero, so ${}^tyMx=x_1w+x_2\alpha(w)$ is nonzero since $V=W_1\oplus W_2$ and $\alpha$ is injective.
	The entries of $Mx$ lie in $W_1$, $W_2$, and $W_3$ for $x={}^t(1,0)$, ${}^t(0,1)$, and ${}^t(1,1)$, respectively.
	This proves the assertion in the case \cref{item:conic-graph}.
\end{proof}

\subsection{Characterization of Fano blowups}\label{ss:grass-fano}

In this subsection, we prove that, for $n\ge3$, the variety $\tilde X$ is Fano if and only if $\theta<n+1$.

We first prove a numerical lemma, which gives the sets of indices needed to apply \cref{lem:grass-line,lem:grass-conic}.

\begin{lem}\label{lem:grass-numerical}
	In \cref{set:grass}, the following hold.
	\begin{enumerate}
		\item\label{item:num-select} Assume that $n\ge3$, and let $I'$ be a set of indices with $\sum_{i\in I'}(r_i-1)\ge n+1$.
			Then there is a subset $I\subseteq I'$ with
			\[
				\sum_{i\in I}(r_i-1)\ge n+1,\qquad\sum_{i\in I}(c_i-1)\le2(n-d).
			\]
		\item\label{item:num-reduce} Assume that $\sum_{i=1}^s(r_i-1)\ge n+1$ and that $\sum_{i\in I}(r_i-1)\le n$ for every set $I$ of indices with $\dim W_I\ge d-1$.
			Then $d=2$, and one of the following holds.
			\begin{itemize}
				\item $r_i=2$ for every $i$.
				\item $r_i=n+1$ for every $i$, and $W_i\cap W_j=0$ for $i\ne j$.
			\end{itemize}
	\end{enumerate}
\end{lem}

\begin{proof}
	We first prove \cref{item:num-select}.
	Renumbering the indices, we may assume that $I'=\{1,\dotsc,t\}$ and $r_1\ge\dotsb\ge r_t$.
	Let $m$ be the smallest integer with $\sum_{i\le m}(r_i-1)\ge n+1$, and set $I\coloneq\{1,\dotsc,m\}$.
	It remains to show that $\sum_{i\in I}(c_i-1)\le2(n-d)$.
	If $m\le2$, then it follows from $c_i\le n+1-d$.
	Assume that $m\ge3$.
	By the minimality of $m$, we have $\sum_{i<m}(r_i-1)\le n$.
	If $d=1$, then $r_i=c_i$, and the ordering gives $2(c_m-1)\le(c_1-1)+(c_2-1)\le n$, so $\sum_{i\in I}(c_i-1)\le n+\lfloor n/2\rfloor\le2(n-1)$, since $n\ge3$.
	Assume that $d\ge2$.
	Then
	\[
		d\sum_{i<m}(c_i-1)=\sum_{i<m}((r_i-1)-(d-1))\le n-(m-1)(d-1)\le d(n-d),
	\]
	where the last inequality is equivalent to $(d-1)(n-d+m-2)\ge1$.
	Together with $c_m-1\le n-d$, this gives $\sum_{i\in I}(c_i-1)\le2(n-d)$.
	This proves \cref{item:num-select}.

	We next prove \cref{item:num-reduce}.
	If $d=1$, then $I=\{1,\dotsc,s\}$ satisfies $\dim W_I\ge d-1$, contrary to the assumption.
	Hence $d\ge2$.
	We first show that every $I$ with $\sum_{i\in I}(r_i-1)\ge n+1$ satisfies
	\begin{equation}\label{eq:num-lower}
		\sum_{i\in I}c_i\ge n+3-d.
	\end{equation}
	Indeed, if $\sum_{i\in I}c_i\le n+2-d$, then $W_I$ has dimension at least $d-1$, contradicting the assumption.
	Since $c_i\le n+1-d$, applying \eqref{eq:num-lower} to $I=\{i\}$ gives $r_i-1\le n$ for every $i$.

	Order the indices so that $r_1\ge\dotsb\ge r_s$.
	Let $m$ be the smallest integer with $\sum_{i\le m}(r_i-1)\ge n+1$, so that $d\le r_m\le n+1$.
	The ordering and the minimality of $m$ give
	\begin{equation}\label{eq:num-min}
		(m-1)(r_m-1)\le\sum_{i<m}(r_i-1)\le n.
	\end{equation}
	Moreover, \eqref{eq:num-lower} gives $d(n+3-d)-m\le\sum_{i\le m}(r_i-1)\le n+r_m-1$.
	Hence, setting $N\coloneq d(n+3-d)-n$, we have $m-1\ge N-r_m$, and therefore $(N-r_m)(r_m-1)\le n$.
	Since $N=(d-2)(n+1-d)+n+2$, this inequality can be rewritten as
	\[
		(d-2)(n+1-d)(r_m-1)+(r_m-2)(n+1-r_m)\le0.
	\]
	Since $2\le d\le r_m\le n+1$, both terms on the left-hand side are nonnegative, so both vanish.
	The vanishing of the first term gives $d=2$.
	The vanishing of the second term gives $r_m\in\{2,n+1\}$.
	Moreover, the inequalities in \eqref{eq:num-min} are equalities, so $r_1=\dotsb=r_m\in\{2,n+1\}$.

	If $r_1=2$, then $r_i=2$ for every $i$, since $2=d\le r_i\le r_1$.
	Assume that $r_1=n+1$.
	If $r_i<n+1=r_1$ for some $i$, then $c_1+c_i<2c_1=n+1$, while $(r_1-1)+(r_i-1)\ge n+1$, which contradicts \eqref{eq:num-lower}.
	Thus $r_i=n+1$ for every $i$.
	For $i\ne j$, we have $(r_i-1)+(r_j-1)=2n>n$, so the assumption with $d=2$ gives $W_i\cap W_j=0$.
	This proves \cref{item:num-reduce}.
\end{proof}

We now prove the criterion, using \cref{prop:grass-positive} and the rational curves of \cref{ss:grass-curves}.

\begin{thm}\label{thm:grass}
	In \cref{set:grass}, assume that $n\ge3$.
	Then $\tilde X$ is Fano if and only if $\theta<n+1$.
\end{thm}

\begin{proof}
	If $\theta<n+1$, then $\tilde X$ is Fano by \cref{prop:grass-positive}.

	Conversely, assume that $\theta\ge n+1$.
	We show that $\tilde X$ is not Fano.
	Assume first that there is a set $I'$ of indices with $\dim W_{I'}\ge d-1$ and $\sum_{i\in I'}(r_i-1)\ge n+1$.
	By \cref{lem:grass-numerical}~\cref{item:num-select}, there is a subset $I\subseteq I'$ with $\sum_{i\in I}(r_i-1)\ge n+1$ and $\sum_{i\in I}(c_i-1)\le2(n-d)$.
	Since $W_{I'}\subseteq W_I$, we have $\dim W_I\ge d-1$, so \cref{lem:grass-line,lem:grass-lift} show that $\tilde X$ is not Fano.

	Assume next that no such set $I'$ exists.
	Since every term in \cref{defi:theta} is at most $\sum_i(r_i-1)$, we have $\sum_i(r_i-1)\ge\theta\ge n+1$.
	By \cref{lem:grass-numerical}~\cref{item:num-reduce}, we have $d=2$, and either $r_i=2$ for every $i$, or $r_i=n+1$ for every $i$ and $W_i\cap W_j=0$ for $i\ne j$.
	Since $d=2$, the terms in \cref{defi:theta} are $\sum_{i\in I}(r_i-1)\le n$ for $W_I\ne0$ and $\sum_{i\in I}(r_i-1)/2$ for $W_I=0$, so $\sum_i(r_i-1)/2\ge\theta\ge n+1$.

	If $r_i=2$ for every $i$, then $s\ge2n+2$.
	Let $I$ be a set of $2n+2$ indices.
	Every subset $I_0\subseteq I$ with $W_{I_0}\ne0$ satisfies $\card{I_0}\le n<\card I/2$, and $\card I\le3n-1$ since $n\ge3$.
	Thus \cref{lem:grass-conic}~\cref{item:conic-quadric} and \cref{lem:grass-lift} show that $\tilde X$ is not Fano.

	If $r_i=n+1$ for every $i$, then $sn\ge2n+2$ gives $s\ge3$.
	By \cref{lem:grass-conic}~\cref{item:conic-graph} with $I=\{1,2,3\}$ and \cref{lem:grass-lift}, there is a morphism $\tilde g\colon\PP^1\to\tilde X$ with $\deg\tilde g^*(-K_{\tilde X})\le2(n+1)-3n<0$.
	Hence $\tilde X$ is not Fano.
\end{proof}

Finally, we make \cref{thm:grass} explicit for each $d$.

\begin{cor}\label{cor:grass}
	In \cref{set:grass}, assume that $n\ge3$.
	Then the following hold.
	\begin{enumerate}
		\item\label{item:grass-dne2} Assume that $d\ne2$.
			Then $\tilde X$ is Fano if and only if $\sum_i(r_i-1)\le n$.
		\item\label{item:grass-d2} Assume that $d=2$.
			Then $\tilde X$ is Fano if and only if $\sum_i(r_i-1)\le2n+1$ and $\sum_{i\in I}(r_i-1)\le n$ for every set $I$ of indices with $W_I\ne0$.
	\end{enumerate}
\end{cor}

\begin{proof}
	For \cref{item:grass-dne2}, if $\sum_i(r_i-1)\le n$, then $\theta\le\sum_i(r_i-1)<n+1$.
	Conversely, assume that $\sum_i(r_i-1)\ge n+1$.
	Since $d\ne2$, \cref{lem:grass-numerical}~\cref{item:num-reduce} shows that there is a set $I$ of indices with $\dim W_I\ge d-1$ and $\sum_{i\in I}(r_i-1)\ge n+1$, and hence $\theta\ge n+1$.
	Thus \cref{item:grass-dne2} follows from \cref{thm:grass}.

	For \cref{item:grass-d2}, the denominator in \cref{defi:theta} is $1$ for $W_I\ne0$ and $2$ for $W_I=0$.
	Since the term for $I=\{1,\dotsc,s\}$ is at least $\sum_i(r_i-1)/2$, we obtain
	\[
		\theta=\max\biggl\{\frac12\sum_{i=1}^s(r_i-1),\ \max_I\sum_{i\in I}(r_i-1)\biggr\},
	\]
	where $I$ runs over the sets of indices with $W_I\ne0$.
	Since $\theta$ is a half-integer, \cref{item:grass-d2} follows from \cref{thm:grass}.
\end{proof}

\begin{rem}\label{rem:grass-n2}
	The assumption $n\ge3$ in \cref{thm:grass} is needed.
	Let $n=2$ and $d=1$, and let $W_1,W_2,W_3$ be one-dimensional subspaces that span $V$.
	Then $\tilde X$ is the blowup of $\Gr(1,V)=\PP^2$ at three points not on a line, which is a del Pezzo surface, while $\theta=\sum_i(r_i-1)=3=n+1$.
\end{rem}

\subsection{Grassmannians of dimension four}\label{ss:grass-fourfolds}

In this subsection, we apply \cref{cor:grass} to the Grassmannians of dimension four, which are $\PP^4=\Gr(1,k^{\oplus5})$ and the smooth quadric fourfold $\Gr(2,k^{\oplus4})$ up to isomorphism.

We first consider planes in $\PP^4$.
They are the subvarieties $\Gr(1,W)$ for three-dimensional subspaces $W\subseteq k^{\oplus5}$.
By \cref{lem:subgrass}, planes in $\PP^4$ meet transversely if and only if any two of them meet in a point and no three of them have a common point.

\begin{cor}\label{cor:P4}
	In \cref{set:grass}, assume that $V=k^{\oplus5}$, $d=1$, and $r_i=2$ for every $i$, so that $\tilde X$ is the successive blowup of $X=\PP^4$ along planes $Z_1,\dotsc,Z_s$.
	Then $\tilde X$ is Fano if and only if $s\le4$.
\end{cor}

\begin{proof}
	Here $n=4$ and $\sum_i(r_i-1)=s$, so the result follows from \cref{cor:grass}~\cref{item:grass-dne2}.
\end{proof}

\Cref{cor:grass} also applies to centers of different dimensions.
In \cref{set:grass} with $V=k^{\oplus5}$ and $d=1$, the centers $\Gr(1,W_i)\subseteq\PP^4$ are points, lines, and planes, since $r_i\ge2$, and they have $r_i-1=3,2,1$.
Hence, if $s\ge3$ and not all the centers are planes, then $\tilde X$ is Fano only for two planes and a line, which we now describe.

\begin{exa}\label{exa:P4-line}
	Let $Z_1,Z_2\subseteq\PP^4$ be planes that meet in a point, and let $Z_3\subseteq\PP^4$ be a line disjoint from $Z_1\cup Z_2$.
	By \cref{lem:subgrass}~\cref{item:subgrass-transverse}, the centers $Z_1,Z_2,Z_3$ meet transversely.
	Since $\sum_i(r_i-1)=1+1+2=4$, the successive blowup $\tilde X$ of $\PP^4$ along them is Fano by \cref{cor:grass}~\cref{item:grass-dne2}, and $\rho(\tilde X)=4$.

	The anticanonical degree of $\tilde X$ is $(-K_{\tilde X})^4=293$, which we compute as follows.
	Let $\tilde Y\to Y$ be the blowup of a smooth projective fourfold along a smooth subvariety $Z$ of codimension $r\in\{2,3\}$ with normal sheaf $N$, and set $A\coloneq\restr{(-K_Y)}{Z}$.
	Computing on the exceptional divisor as in the proof of \cite[Lemma~3.2]{CasagrandeRomano}, we get
	\[
		(-K_{\tilde Y})^4=
		\begin{cases}
			(-K_Y)^4-6A^2+4A\cdot c_1(N)-c_1(N)^2+c_2(N) & (r=2),\\
			(-K_Y)^4-32\deg A+16\deg N & (r=3).
		\end{cases}
	\]
	First, $(-K_{\PP^4})^4=625$.
	On $Z_1\isom\PP^2$, we have $A=5h$ and $N\isom\Oc(1)^{\oplus2}$ for the class $h$ of a line, so blowing up $Z_1$ lowers the degree by $113$.
	Since $Z_1$ and $Z_2$ meet transversely, the strict transform of $Z_2$ coincides with its inverse image.
	Hence it is the blowup of $Z_2$ at the point $Z_1\cap Z_2$, and its normal sheaf is the pullback of $N_{Z_2/\PP^4}\isom\Oc(1)^{\oplus2}$.
	On it, $A=5h-e$ for the pullback $h$ of the class of a line and the exceptional curve $e$, so blowing it up lowers the degree by $107$.
	The line $Z_3$ is disjoint from the first two centers, and $A$ has degree $5$ and $N\isom\Oc(1)^{\oplus3}$ on it, so blowing it up lowers the degree by $112$.
	Hence
	\[
		(-K_{\tilde X})^4=625-113-107-112=293.
	\]
	In particular, $\tilde X$ is not a product of surfaces, since a product of del Pezzo surfaces of degrees $d_1,d_2$ has the even anticanonical degree $6d_1d_2$.

	Over $\CC$, the blowup formula for cohomology shows that $h^{p,q}(\tilde X)=0$ for $p\ne q$, that $h^{1,1}(\tilde X)=h^{3,3}(\tilde X)=4$, and that $h^{2,2}(\tilde X)=6$.
\end{exa}

We next consider planes in a smooth quadric fourfold $X=\Gr(2,V)\subseteq\PP((\bigwedge^2V)^\vee)=\PP^5$, where $\dim V=4$.
This covers all smooth quadric fourfolds up to projective equivalence.
The planes in $X$ form two families.
One consists of the planes $\Gr(2,W)$ for hyperplanes $W\subseteq V$, and the other consists of the planes of points $[P]$ with $P$ containing a fixed one-dimensional subspace of $V$.
The isomorphism $\Gr(2,V)\isomarrow\Gr(2,V^\vee)$ given by $P\mapsto P^\perp$ exchanges the two families.
Hence, up to this isomorphism, planes of one family are the planes $\Gr(2,W)$ for hyperplanes $W\subseteq V$, and such a plane corresponds to the point $[W^\perp]$ of $\PP^3\coloneq\PP(V)$.

By \cref{lem:subgrass}, these planes meet transversely if and only if the corresponding points in $\PP^3$ are pairwise distinct and no three of them are collinear.

\begin{cor}\label{cor:quadric}
	In \cref{set:grass}, assume that $n=3$, $d=2$, and $r_i=2$ for every $i$, so that the $W_i$ are hyperplanes and $\tilde X$ is the successive blowup of the smooth quadric fourfold $X$ along planes $Z_1,\dotsc,Z_s$ of one family.
	Set $p_i\coloneq[W_i^\perp]\in\PP^3=\PP(V)$.
	Then $\tilde X$ is Fano if and only if $s\le7$ and no four of the points $p_i$ lie on a plane.
\end{cor}

\begin{proof}
	We have $\sum_i(r_i-1)=s$ and $\sum_{i\in I}(r_i-1)=\card I$.
	For a set $I$ of indices, $W_I\ne0$ if and only if the one-dimensional subspaces $W_i^\perp\subseteq V^\vee$ with $i\in I$ do not span $V^\vee$, that is, the points $p_i$ with $i\in I$ lie on a plane.
	Thus the result follows from \cref{cor:grass}~\cref{item:grass-d2}.
\end{proof}

\printbibliography

\end{document}